\documentclass[a4j,12pt]{article}
\title{{\sc New Proofs of Triangle Inequalities}}
\author{Norihiro Someyama\footnote{Shin-yo-ji Temple, 5-44-4 Minamisenju, Arakawa-ku, Tokyo 116-0003 Japan / e-mail: {\tt philomatics@outlook.jp}}\quad $\&$\ Mark Lyndon Adamas Borongan\footnote{Mathematics Department, Silliman University, 1 Hibbard Avenue, Dumaguete City, Negros Oriental, 6200 Philippines / e-mail: {\tt lyndonaborongan@su.edu.ph}}}
\date{\empty \vspace{-12mm}}
\usepackage{latexsym,amsmath,amssymb,amsthm,amscd,mathrsfs}
\usepackage{bm}
\usepackage[pdftex]{graphicx}
\usepackage{fancybox}
\usepackage{ascmac}
\usepackage{multicol}
\usepackage{tabularx}
\usepackage{comment}
\usepackage[top=4cm,bottom=4cm,left=3cm,right=3cm]{geometry}
\theoremstyle{theorem}
\newtheorem{thm}{Theorem}[section]
\newtheorem{prop}[thm]{Proposition}

\newtheorem{lem}[thm]{Lemma}

\theoremstyle{definition}

\newtheorem{assump}{Assumption}[section]
\newtheorem{rem}{Remark}[section]

\theoremstyle{theorem}

\makeatletter
 
 \@addtoreset{equation}{section}
\makeatother

\begin{document}
\maketitle

\begin{abstract}
We give three new proofs of the triangle inequality in Euclidean Geometry. 
There seems to be only one known proof at the moment.
It is due to properties of triangles, but our proofs are due to circles or ellipses.
We aim to prove the triangle inequality as simple as possible without using properties of triangles.
\end{abstract}

{\small
{\bf Keywords and Phrases}: triangle inequality, circle, ellipse, focci of ellipse
}

\section{Introduction}
We consider triangle inequalities for triangles in Euclidean space.
For convenience, let the dimension of the space be $2$.
We write $\triangle{\rm ABC}$ for the triangle with three vertices ${\rm A},{\rm B},{\rm C}\in \mathbb R^2$.
Hereafter, ${\rm AB}$ denotes the segment from a point ${\rm A}$ to a point ${\rm B}$, and $\angle{\rm ABC}$ the angle made by segments ${\rm AB},{\rm BC}$.
Moreover, $\overline{{\rm PQ}}$ stands for the length of the segment ${\rm PQ}$ in $\mathbb R^2$.

The triangle inequality asserts that the sum of any two sides of a triangle is {\it strictly} bigger than the remaining third side.
This geometric inequality is well known as one of the most fundamental and classical theorems in Euclidean geometry:

\begin{thm}[{\bf Triangle Inequalities}]
\label{thm:CTI}
For any triangle $\triangle{\rm ABC}$, an inequality
\begin{equation}
\overline{{\rm AB}}+\overline{{\rm AC}}>\overline{{\rm BC}}
\label{eq:CTI}
\end{equation}
holds (regardless of the dimension of the space).
\end{thm}

This was probably proved by the ancient Greeks for the first time, but the proof is still considered important.
See Subsection \ref{subsec:EG} for details of `the first proof'.
However, we want to give new and more natural proofs in this article.
Furthermore, it is also well known that we can prove the triangle inequality {\it in the broad sence}, i.e. the `$\ge $'-version of (\ref{eq:CTI}), by algebraic argument.
It is not exactly the triangle inequality in the sense of Euclidean geometry, because the point ${\rm A}$ is on the segment ${\rm BC}$ in case $\overline{{\rm AB}}+\overline{{\rm AC}}=\overline{{\rm BC}}$.
For details of them, Subsection \ref{subsec:ATI} (in particular Remark \ref{rem:TIrem}) will mention.

Throughout this article, we always assume the following:

\begin{assump}
\label{assump:BC>ABAC}
${\rm BC}$ is the longest side among three sides of $\triangle {\rm ABC}$. 
That is,
\begin{equation}
\label{eq:BC>ABAC}
\overline{{\rm BC}}>\max\{\overline{{\rm AB}},\overline{{\rm AC}}\}.
\end{equation}
\end{assump}

\subsection{The best known Proof by Euclidean Geometrical Method\label{subsec:EG}}
\begin{figure}[h]
\begin{center}
\includegraphics[width=7cm]{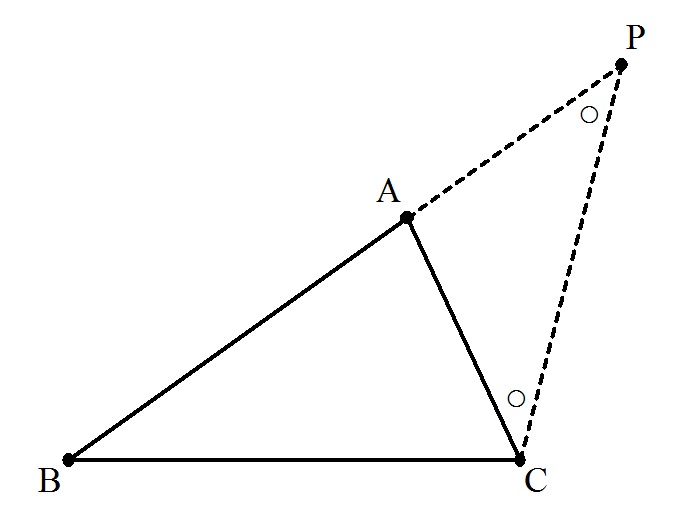}
\end{center}
\caption{}
\label{fig:ABC}
\end{figure}

We consider a triangle $\triangle {\rm ABC}$ with Assumption \ref{assump:BC>ABAC} and plot a point ${\rm P}$, obeying 
\begin{align}
\label{eq:EG1}
\overline{{\rm AC}}=\overline{{\rm AP}}, 
\end{align}
on an extension of the segment ${\rm BA}$.
(See Figure \ref{fig:ABC}.)
Then, 
\begin{align}
\label{eq:EG2}
\angle{\rm BCP}>\angle{\rm ACP}=\angle{\rm APC}=\angle{\rm BPC}.
\end{align}
Hence, $\triangle{\rm PBC}$ satisfies
\[
\overline{{\rm BC}}<\overline{{\rm BP}}=\overline{{\rm BA}}+\overline{{\rm AP}}=\overline{{\rm AB}}+\overline{{\rm AC}},
\]
from (\ref{eq:EG1}) and (\ref{eq:EG2}).
This completes the proof.

{\small
\begin{rem}
The above proof also applies to obtuse angled triangles, but there is another short proof for obtuse angled triangles.
In fact, if $\angle{\rm BCA}$ is the maximum angle, 
\[
\overline{{\rm AB}}>\max\{\overline{{\rm BC}},\overline{{\rm AC}}\}
\]
holds, so we obviously gain (\ref{eq:CTI}).
\end{rem}
}

\subsection{Proofs of Algebraic Triangle Inequalities \label{subsec:ATI}}
The names `triangle inequalities' often appear in some fields such as linear algebra and functional analysis.
These `triangle inequalities' are able to argue by algebraic method.

\subsubsection{Triangle Inequalities for Absolute Values of Vectors \label{subsubsec:AG}}
We consider three vectors 
\[
{\bm a}:=\overrightarrow{{\rm BA}}, \quad
{\bm b}:=\overrightarrow{{\rm AC}},\quad
{\bm c}:=\overrightarrow{{\rm BC}}.
\]
Then, since ${\bm c}={\bm a}+{\bm b}$, it is sufficient to prove
\begin{equation}
|{\bm a}+{\bm b}|<|{\bm a}|+|{\bm b}|
\label{eq:abab}
\end{equation}
so as to obtain (\ref{eq:CTI}), where $|\quad |$ is an Euclidean norm:
\[
|{\bm x}|:=\sqrt{x_1^2+x_2^2}\qquad {\rm for}\quad  {\bm x}=(x_1,x_2).
\]
Note that, of course, $|{\bm a}|,|{\bm b}|>0$.
To see (\ref{eq:abab}), we consider the following equation:
\begin{align}
\label{eq:ab-ab}
(|{\bm a}|+|{\bm b}|)^2-|{\bm a}+{\bm b}|^2=2(|{\bm a}||{\bm b}|-{\bm a}\cdot {\bm b})
\end{align}
where `$\,\cdot\,$' is the inner product.
The inner product defines by
\[
{\bm a}\cdot {\bm b}
=|{\bm a}||{\bm b}|\cos(\pi-\angle{\rm BAC})
=-|{\bm a}||{\bm b}|\cos \angle{\rm BAC}
\]
where $0^{\circ}<\angle{\rm BAC}<180^{\circ}$.
Thus, 
\begin{align}
\label{eq:ab-ab>0}
|{\bm a}||{\bm b}|-{\bm a}\cdot {\bm b}
=|{\bm a}||{\bm b}|(1+\cos\theta)>0,
\end{align}
since $0\le \cos\theta<1$.
Hence, (\ref{eq:abab}) is proved from (\ref{eq:ab-ab}) and (\ref{eq:ab-ab>0}).

{\small
\begin{rem}
It follows in more general that
\begin{align}
|{\bm a}+{\bm b}|\le |{\bm a}|+|{\bm b}|
\label{eq:abab'}.
\end{align}
Let us call this the {\it algebraic triangle inequality}.
To see (\ref{eq:abab'}), we consider (\ref{eq:ab-ab}) and the Cauchy-Schwarz-Bunyakovsky's inequality
\begin{align*}
|{\bm a}\cdot {\bm b}|\le |{\bm a}||{\bm b}|.
\end{align*}
They immediately derive (\ref{eq:abab'}).
\end{rem}
}

\subsubsection{Algebraic Triangle Inequalities for Complex Numbers \label{subsubsec:C}}
(\ref{eq:abab}) can be shown by using complex values.
We write 
\[
{\bm a}:=(a_1,a_2),\quad \ {\bm b}:=(b_1,b_2)
\]
for two real vectors in Subsection \ref{subsubsec:AG}.
Since the mapping ${\mathbb R}^2\ni (x,y)\mapsto x+iy\in \mathbb C$ is isomorphic: $\mathbb R^2\approx \mathbb C$, let ${\bm a},{\bm b}$ correspond to two complex numbers
\[
z_1:=a_1+ia_2=r_1e^{i\theta_1},\quad \ z_2:=b_1+ib_2=r_2e^{i\theta_2}
\]
respectively.
Here $i:=\sqrt{-1}$ and $0^{\circ}<\theta_j<180^{\circ}\ (j=1,2)$.
Note that $\theta_1\neq \theta_2$.
Then, we should prove
\begin{equation}
\label{eq:z1z2}
|z_1+z_2|<|z_1|+|z_2|
\end{equation}
to obtain (\ref{eq:abab}).
We can see (\ref{eq:z1z2}) as follows:
\begin{align*}
|z_1+z_2|&=|(r_1\cos \theta_1+r_2\cos \theta_2)+i(r_1\sin \theta_1+r_2\sin \theta_2)| \\
&=\sqrt{(r_1\cos \theta_1+r_2\cos \theta_2)^2+(r_1\sin \theta_1+r_2\sin \theta_2)^2} \\
&=\sqrt{r_1^2+r_2^2+2r_1r_2\cos(\theta_1-\theta_2)} \\
&<\sqrt{r_1^2+2r_1r_2+r_2^2} \\
&=r_1+r_2 \\
&=|z_1|+|z_2|.
\end{align*}
We here use that $|\cos(\theta_1-\theta_2)|<1$.
Hence, this completes the proof.

{\small
\begin{rem}
\label{rem:TIrem}
More generally, we have 
\begin{align}
\label{eq:z1z2'}
|z_1+z_2|\le |z_1|+|z_2|.
\end{align}
This corresponds to (\ref{eq:abab'}).
The condition for (\ref{eq:abab'}) (resp. (\ref{eq:z1z2'})) to become an equation is what ${\bm a}$ and ${\bm b}$ are parallel (resp. that $\theta:=\theta_1-\theta_2=n\pi$ for any $n\in \mathbb Z$).
Geometrically, triangles are of course not made if these conditions are satisfied, so the algebraic triangle inequalities such as (\ref{eq:abab}) and $(\ref{eq:z1z2})$ are triangle inequalities in the broad sence.
\end{rem}
}

\section{New Proofs by Circles\label{sec:prfc}}
Let us present our proofs of triangle inequalities from this section.
This section gives us two proofs of (\ref{eq:CTI}) by circles.

\subsection{New Proof 1 by Circles\label{subsec:Prf1c}}
\begin{figure}[h]
\begin{center}
\includegraphics[width=8.5cm]{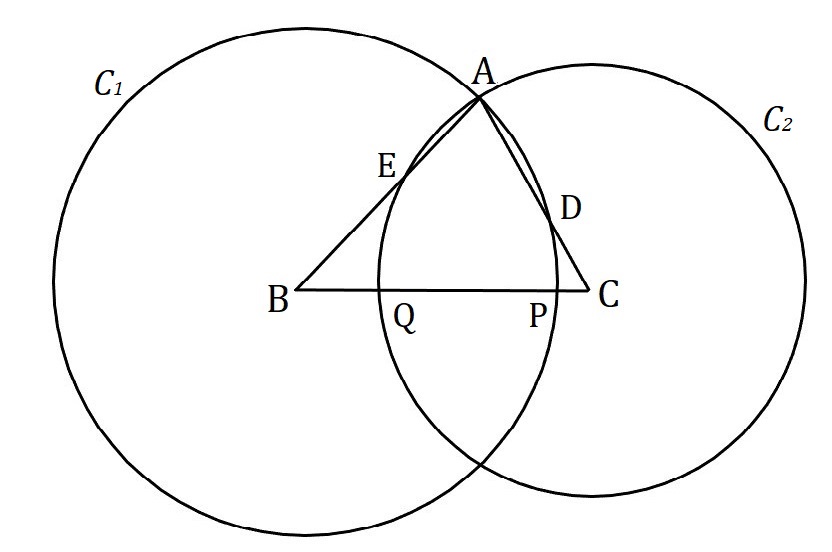}
\end{center}
\caption{}
\label{fig:ABCcir}
\end{figure}

We consider two circles, one is a circle $C_1$ whose center is the point ${\rm B}$ and whose radius is the segment ${\rm AB}$, and the other is a circle $C_2$ whose center is the point ${\rm C}$ and whose radius is the segment ${\rm AC}$.
(See Figure \ref{fig:ABCcir}.)
If 
\[
\overline{{\rm AB}}+\overline{{\rm AC}}\le \overline{{\rm BC}},
\]
then $C_1$ and $C_2$ come in contact with each other or become separated from each other.
In either case, segments ${\rm AB}$, ${\rm BC}$ and ${\rm AC}$ never make a triangle.
Hence, (\ref{eq:CTI}) must hold.

\subsection{New Proof 2 by Circles \label{subsec:Prf2c}}
If $C_1$ and $C_2$ come in contact with each other, the point ${\rm A}$ is on the segment ${\rm BC}$.
So, it is sufficient to argue only the case that $C_1$ intersects $C_2$ at two points.
Then, ${\rm A}$ is one of two points of intersection of $C_1$ and $C_2$.
(See Figure \ref{fig:ABCcir}.)
Remark that the point ${\rm C}$ (resp. ${\rm B}$) is outside $C_1$ (resp. $C_2$) from (\ref{eq:BC>ABAC}).
We write ${\rm P}$ for a point of intersection of ${\rm BC}$ and $C_1$.
Similarly, we write ${\rm Q}$ for a point of intersection of ${\rm BC}$ and $C_2$.
Then,
\[
\overline{{\rm BC}}=\overline{{\rm BP}}+\overline{{\rm PC}}=\overline{{\rm AB}}+\overline{{\rm PC}}
\]
because $\overline{{\rm AB}}=\overline{{\rm BP}}$.
So it is sufficient to prove
\begin{align}
\label{eq:AC>PC}
\overline{{\rm AC}}>\overline{{\rm PC}},
\end{align}
but, since $C_1$ intersects $C_2$, we obviously obtain (\ref{eq:AC>PC}) as follows:
\[
\overline{{\rm PC}}<\overline{{\rm QC}}=\overline{{\rm AC}}.
\]
Hence, this completes the proof.

{\small
\begin{rem}
A proof for obtuse angled triangles is known and is similar to the above proof. 
(See \cite{D} for example.)
That is, if $\angle{\rm BAC}$ is the maximum angle, we consider a circle whose center is ${\rm C}$ and whose radius is ${\rm BC}$.
We leave the details to the reader.
\end{rem}
}

\section{New Proof by Ellipses and Properties of $\mathbb R$\label{sec:ellip}}
We first give a {\it partial} proof of the triangle inequality by ellipses in Section \ref{subsec:ge2a}.
We next complement the proof by certain property of real numbers in Section \ref{subsec:<2a}.
That complementary proof is however only an appendix.

\subsection{Proof for Special Triangles \label{subsec:ge2a}}
\begin{figure}
\begin{center}
\includegraphics[width=8.5cm]{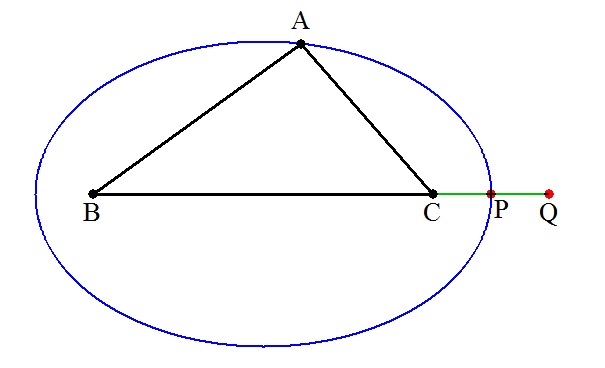}
\end{center}
\caption{}
\label{fig:Elli}
\end{figure}

We consider an ellipse
\[
{\cal E}:\ \frac{x^2}{a^2}+\frac{y^2}{b^2}=1\quad (a>b>0)
\]
and set that vertices ${\rm B}$ and ${\rm C}$ of $\triangle{\rm ABC}$ are focci of ${\cal E}$ given by 
\[
{\rm B}(-\sqrt{a^2-b^2},0),\quad {\rm C}(\sqrt{a^2-b^2},0).
\]
Then, if the vertex ${\rm A}$ of $\triangle{\rm ABC}$ is any point on ${\cal E}$, we gain
\begin{equation}
\overline{{\rm AB}}+\overline{{\rm AC}}=2a
\label{eq:=2a}
\end{equation}
by definition of an ellipse.
If ${\rm P}(a,0)$, it satisfies that 
\begin{equation}
\overline{{\rm AB}}+\overline{{\rm AC}}=\overline{{\rm BP}}+\overline{{\rm PC}}\ (=2a)
\label{eq:ABACPC}
\end{equation}
from (\ref{eq:=2a}).
We plot a point ${\rm Q}$, on an extension of the segment ${\rm BP}$, obeying 
\begin{align}
\label{eq:PCPQ}
\overline{{\rm PC}}=\overline{{\rm PQ}}.
\end{align}
(See Figure \ref{fig:Elli}.) 
Hence, by (\ref{eq:ABACPC}) and (\ref{eq:PCPQ}),
\[
\overline{{\rm AB}}+\overline{{\rm AC}}=\overline{{\rm BP}}+\overline{{\rm PQ}}=\overline{{\rm BQ}}>\overline{{\rm BC}},
\]
so (\ref{eq:CTI}) holds if (\ref{eq:=2a}).

From the above, we can also see that (\ref{eq:CTI}) holds if
\begin{equation}
\overline{{\rm AB}}+\overline{{\rm AC}}>2a.
\label{eq:>2a}
\end{equation}

Thus, we obtain the following result:

\begin{prop}
\label{prop:ge2a}
Let $a>0$ be the semi-major axis of ${\cal E}$.
If 
\[
\overline{{\rm AB}}+\overline{{\rm AC}}\ge 2a,
\]
then (\ref{eq:CTI}) holds.
\end{prop}

{\small
\begin{rem}
\begin{itemize}
\item[1)] The above proof of Proposition \ref{prop:ge2a} is the argument by analytic geometry, but we can also prove by the algebraic argument as follows:
\[
\overline{{\rm AB}}+\overline{{\rm AC}}\ge 2a>2\sqrt{a^2-b^2}=\overline{{\rm BC}}.
\]
\item[2)] The above proof is acceptable wheather $\triangle{\rm ABC}$ is acute or not.
\end{itemize}
\end{rem}
}

\subsection{Appendix: Completion of the Proof \label{subsec:<2a}}
As we saw in Section \ref{subsec:ge2a}, it was easy to see (\ref{eq:CTI}) if
\[
\overline{{\rm AB}}+\overline{{\rm AC}}\ge 2a.
\]
This is a partial proof of (\ref{eq:CTI}).
In order to complete the proof of (\ref{eq:CTI}), it is necessary to prove the following statement.

\begin{prop}
Let $a>0$ be the semi-major axis of ${\cal E}$.
If 
\[
\overline{{\rm AB}}+\overline{{\rm AC}}<2a,
\]
then (\ref{eq:CTI}) holds.
\label{prop:<2a}
\end{prop}

To see this, we use the following property for real numbers.

\begin{lem}[e.g. \cite{B}, p.42]
\label{lem:acbcab}
Let $\alpha,\beta,\gamma\in \mathbb R$.
If $\alpha<\gamma$ for any $\gamma$ such that $\beta<\gamma$, then $\alpha\le \beta$.
\end{lem}

{\small
\begin{rem}
Lemma \ref{lem:acbcab} is a statement which paraphrases the following well-known proposition used frequently in measure theory: 
\begin{center}
`{\it Let $a,b\in \mathbb R$. 
If $a<b+\varepsilon$ for any $\varepsilon>0$, then $a\le b$.}'
\end{center}
To see Lemma \ref{lem:acbcab} directly, we use reductio ad absurdum. 
In fact, we should assume $\alpha>\beta$ and set $\gamma=(\alpha+\beta)/2$.
\end{rem}
}

Actually, we can immediately prove Proposition \ref{prop:<2a} by putting
\[
\alpha=\overline{{\rm BC}},\quad \ 
\beta=\overline{{\rm AB}}+\overline{{\rm AC}},\quad \ 
\gamma=2a
\]
in Lemma \ref{lem:acbcab}.

Hence, by virtue of Proposition \ref{prop:ge2a} and Proposition \ref{prop:<2a}, this completes the proof of (\ref{eq:CTI}).
\vspace{4mm}

{\small
\begin{center}
{\bf Acknowledgement}
\end{center}
The authors NS and MLAB thank Honorary Prof. Shigeru Nakamura of Tokyo University of Marine Science and Technology, Dr. Millard R. Mamhot of Silliman University and other doctors, for accurate advice and suggestion.
}

\end{document}